\documentclass[12pt]{amsart}
\usepackage{amsmath,amssymb, amsfonts,amscd,
psfrag,graphicx}

\newcommand{\lebn}

\theoremstyle{plain}

\usepackage[utf8]{inputenc} 
\usepackage[T1]{fontenc} 

\theoremstyle{definition}

\begin{document}

\bibliographystyle{plain}

\title[Lie Derivations of a Matrix Ring over an Associative Ring ]{ Lie
Derivations of a Matrix Ring over an Associative Ring }
\author{ UMUT SAYIN }
\address{ Department of Mathematics, Düzce University, 81620 Konuralp,Düzce, Turkey}
\email{umutsayin@duzce.edu.tr}
\author{ FERİDE KUZUCUOĞLU }
\address{ Department of Mathematics, Hacettepe University, 06800
Beytepe,Ankara, Turkey }
\email{feridek@hacettepe.edu.tr}

\begin{abstract}
Let $K$\ be a 2-torsion free ring with identity. We give a description of
the Lie derivations of $R=R_{n}(K,J)=NT_{n}(K)+M_{n}(J),$\ the ring of all $%
n\times n$\ matrices over $K$\ such that the entries on and above the main
diagonal are elements of an ideal $J$\ of $K.$
\end{abstract}

\maketitle

\section{Introduction }

Let $K$\ be an associative ring with identity. An additive map $\Delta
:K\rightarrow K$\ is called a Lie derivation of $K$\ if $\Delta (x\ast
y)=\Delta (x)\ast y+x\ast \Delta (y)$\ for all $x,y\in K$\ where $x\ast
y=xy-yx.$\ Clearly, derivations are basic examples of Lie derivations of $K.$%
\ The usual problem is to determine whether a Lie derivation is a
derivation. The problem of describing all Jordan derivations, Lie
derivations and Lie automorphisms of a ring or algebra was considered by
many authors(see \cite{BeiBreChe}, \cite{Dragomir}, \cite{Kuzucuoglu}, \cite{LevcOksa}, \cite{Martindale}, \cite{OuWangYao}).
Derivations and Jordan derivations of the ring $R$\ are described in \cite%
{KuzuSayi} and \cite{SayiKuzu}, respectively. In this paper, we continue the
investigation of Lie derivations of $R=R_{n}(K,J)=NT_{n}(K)+M_{n}(J)$\ for $%
n\geq 3$\ where $NT_{n}(K)$\ is the niltriangular $n\times n$\ matrix ring\
and $M_{n}(J)$\ is the ring of all $n\times n$\ matrices over an ideal $J$\
of $K.$\ If $[x_{i,j}]$\ is an element of $R,$\ then $x_{i,j}\in I_{i,j}$\
where $I_{i,j}$\ is equal to $J$\ for $i\leq j$\ or to $R$\ for $i>j.$

The group of automorphisms of this ring $R$\ is given in \cite{KuzuLevc}
under certain restrictions on $J.$\ Throughout this paper, we consider
the square matrices of order $n\geq 3$\ \ and we also assume that the 2-torsion
free ring $K$\ has identity element. Let $e_{i,j}$\ denote the matrix unit
with 1 in the (i,j) position and 0 in every other position. The ring $R$\ is
generated by the sets $Ke_{i+1,i}$\ $(i=1,2,...,n-1)$\ and $Je_{1,n}.$

Let $E$\ denote the identity matrix of order n, $\delta _{i,j}=0$\ for $%
i\neq j$\ and $\delta _{i,i}=1.$\ Then $e_{i,j}\ast e_{r,s}=\delta
_{j,r}e_{i,s}-\delta _{s,i}e_{r,j}.$\ \

\section{Basic Derivations and the Main Theorem}

Let $R=R_{n}(K,J).$\ Then the sets of annihilator of $R$\ and of $J$\ in $K$\ are%
\begin{equation*}
AnnR=\left\{ \alpha \in R:\alpha R=R\alpha =0\right\} 
\end{equation*}%
and%
\begin{equation*}
Ann_{K}(J)=\left\{ \alpha \in K:\alpha J=J\alpha =0\right\},
\end{equation*}%
respectively. By Lemma 1.3 in \cite{KuzuLevc}, we know that $AnnR$\ is equal to $%
Ann_{K}(J)e_{n,1}.$\ Let $C(K)$\ be the center of a ring $K.$\ Clearly, $%
C(R)=AnnR+(J\cap C(K))E.$\ 

\bigskip 

If we replace the original product by the Lie product $x\ast y=xy-yx$\ in
any ring,\ the additive group of that ring\ becomes a Lie ring. If this ring
is of characteristic 2$,$\ then the Lie product coincides with the Jordan
product $x\circ y=xy+yx.$\ So we assume that our ring $K$\ is of
characteristic not 2.

\bigskip 

For any ring $K$\ and any element $x$\ of this ring, the map $\delta
_{x}:K\rightarrow K$\ defined by the rule $\delta _{x}(y)=xy-yx$\ is a
derivation. Such derivations are called inner.

\bigskip 

If $\theta :K\rightarrow K$\ is a derivation of $K$\ satisfying $\theta
(J)\subset J,$\ then\ $\bar{\theta}:[a_{i,j}]\rightarrow \lbrack \theta
(a_{i,j})]$\ is a derivation of $R.$\ These types of derivations are called
ring derivations.

\bigskip 

Let $d=\sum\limits_{i=1}^{n}d_{i}e_{i,i}$\ $(d_{i}\in K).$\ For any
diagonal matrix $d$\ and $x\in R,$\ then the map $\sigma _{d}:x\rightarrow
dx-xd$\ is a derivation called diagonal derivation induced by the matrix $d.$

\bigskip From now on

\begin{itemize}
\item $x_{i,j}$, $y_{i,j}$\ will denote arbitrary elements in $I_{i,j}$\ for
any $i$\ and $j,$

\item $x,$\ $x_{1},$\ $x_{2}$\ will denote arbitrary elements in $K,$

\item $y$, $z$\ will denote arbitrary elements in $J.$
\end{itemize}

Now we construct several types of Lie derivations on $R$\ which will be used
in order to generate\ all Lie derivations$.$

\bigskip 

\textbf{(A)}\ Let $\varsigma _{n}:J\rightarrow Ann_{K}(J)$, \ $\varsigma
_{i}:K\rightarrow Ann_{K}(J),\ \ \alpha _{n}:J\rightarrow C(K),$\ \ $\alpha
_{i}:K\rightarrow C(K)$ and\ $u:J\rightarrow Ann_{K}(J)$\ be additive maps
satisfying%
\begin{equation*}
\begin{tabular}{l}
$(1)\text{ }\varsigma _{i}(J)=\alpha _{i}(J)=0,\text{ }i=1,2,...,n-1,$ \\ 
$(2)\text{ }\varsigma _{n}(J^{2})=\alpha _{n}(J^{2})=0,$ \\ 
$(3)$ $u(xy)=u(yx).$%
\end{tabular}%
\end{equation*}%
Then the map%
\begin{equation*}
\begin{array}{llll}
\Omega : & R & \rightarrow  & C(R) \\ 
& [x_{i,j}] & \rightarrow  & 
\begin{array}{l}
\left( \sum\limits_{i=1}^{n-1}\alpha _{i}(x_{i+1,i})+\alpha
_{n}(x_{1,n})\right) E \\ 
+\left( \sum\limits_{j=1}^{n}u(x_{j,j})+\sum\limits_{i=1}^{n-1}\varsigma
_{i}(x_{i+1,i})+\varsigma _{n}(x_{1,n})\right) e_{n,1}%
\end{array}%
\end{array}%
\end{equation*}%
is a Lie derivation. This Lie derivation is called central since $\Omega =0$%
\ modulo $C(R).$

\bigskip 

We now define special Lie derivations of $R;$

\bigskip 

\textbf{(B1)}\ Special Lie derivation of type I:

The map $%
\begin{array}{llll}
\bar{\Delta}: & xe_{2,1} & \rightarrow  & \gamma (x)e_{n,2}\text{ \ \ \ \ \
\ \ } \\ 
& xe_{n,n-1} & \rightarrow  & \theta (x)e_{n-1,1}\text{ \ \ \ \ \ \ \ } \\ 
& x_{i,j}e_{i,j} & \rightarrow  & 0\text{ \ \ }((i,j)\neq (2,1),(n,n-1))%
\end{array}%
$

is a Lie derivation of $R$\ if the additive maps $\gamma ,\theta
:K\rightarrow Ann_{K}(J)$\ satisfy the conditions%
\begin{equation*}
\begin{array}{l}
i)\text{ }\gamma (J)=\theta (J)=0 \\ 
ii)\text{ }\gamma (x_{1})x_{2}=\gamma (x_{2})x_{1} \\ 
iii)\text{ }x_{1}\theta (x_{2})=x_{2}\theta (x_{1}).%
\end{array}%
\end{equation*}

\textbf{(B2)}\ Special Lie derivation of type II:

Let $\alpha ,\beta ,\gamma :J\rightarrow Ann_{K}(J)$\ be additive maps with
the conditions $\alpha (J^{2})=\beta (J^{2})=\gamma (J^{2})=0,$\ $\alpha
(yx)=x\alpha (y),$\ $\beta (yx)=x\beta (y),$\ $\beta (xy)=\beta (y)x$\ and$\
\gamma (xy)=\gamma (y)x.$\ Then the map%
\begin{equation*}
\begin{array}{llll}
\bar{\Gamma}: & R & \rightarrow  & R \\ 
& ye_{1,n} & \rightarrow  & \alpha (y)e_{n-1,1}+\beta (y)e_{n-1,2}+\gamma
(y)e_{n,2} \\ 
& ye_{1,n-1} & \rightarrow  & -\alpha (y)e_{n,1}-\beta (y)e_{n,2} \\ 
& ye_{2,n} & \rightarrow  & -\beta (y)e_{n-1,1}-\gamma (y)e_{n,1}\text{ ,\ \ 
}y\epsilon J \\ 
& ye_{2,n-1} & \rightarrow  & \beta (y)e_{n,1} \\ 
& x_{i,j}e_{i,j} & \rightarrow  & 0\text{ \ \ \ }((i,j)\neq
(1,n),(1,n-1),(2,n-1),(2,n))%
\end{array}%
\end{equation*}%
is a Lie derivation of $R.$

\bigskip 

\textbf{(B3)}\ Special Lie derivation of type III:

Let $\alpha ,\beta ,\gamma :J\rightarrow Ann_{K}(J)$\ be additive maps
satisfying $\alpha (J^{2})=\gamma (J^{2})=0,$\ $\alpha (yx)=x\alpha (y)$\ and%
$\ \gamma (xy)=\gamma (y)x.$\ Then the map%
\begin{equation*}
\begin{array}{llll}
\Lambda : & R_{3}(K,J) & \rightarrow  & R_{3}(K,J) \\ 
& ye_{1,3} & \rightarrow  & \alpha (y)e_{2,1}+\gamma (y)e_{3,2} \\ 
& ye_{1,2} & \rightarrow  & -\alpha (y)e_{3,1} \\ 
& ye_{2,3} & \rightarrow  & -\gamma (y)e_{3,1} \\ 
& x_{i,j}e_{i,j} & \rightarrow  & 0\text{ \ \ \ }((i,j)\neq
(1,3),(1,2),(2,3))%
\end{array}%
\end{equation*}%
is a Lie derivation of $R_{3}(K,J).$

\bigskip 

\textbf{(C)}\ Let $\alpha ,\beta :J\rightarrow J$\ be additive maps
satisfying the conditions $\alpha (xy)=x\alpha (y),$\ $\beta (yx)=\beta (y)x$%
\ and $\alpha (y)z+y\beta (z)=0.$\ Then the map%
\begin{equation*}
\begin{array}{llll}
\Gamma : & R & \rightarrow  & R \\ 
& ye_{1,n} & \rightarrow  & \alpha (y)e_{1,1}+\beta (y)e_{n,n} \\ 
& ye_{i,n} & \rightarrow  & \alpha (y)e_{i,1}\text{ \ , \ \ }1<i\leq n \\ 
& ye_{1,j} & \rightarrow  & \beta (y)e_{n,j}\text{ \ , \ \ }1\leq j<n \\ 
& x_{i,j} & \rightarrow  & 0\text{ \ \ ,\ \ \ }(i\neq 1\text{, }j\neq n)%
\end{array}%
\end{equation*}%
determines a derivation of the ring $R$\ called almost annihilator
derivation.

\bigskip 

\textbf{Main Theorem:}\ Let $K$\ be a 2-torsion free ring with identity and $%
R=NT_{n}(K)+M_{n}(J)$\ for $n>3.$\ Then every Lie derivation of $R$\ is a
sum of diagonal, inner, ring, almost annihilator derivations and central
Lie, special Lie derivations defined in (B1) and (B2).

\bigskip 

In order to prove our main theorem, we give some technical lemmas. Suppose
that $\Delta $\ is a Lie derivation of $R.$\ For any elementary matrix $%
x_{i,j}e_{i,j}\in R$\ $(x_{i,j}\in I_{i,j},$\ $1\leq i,j\leq n)$, we will
check the properties of $\Delta
(x_{i,j}e_{i,j})=\sum\limits_{s,t=1}^{n}\Delta _{s,t}^{i,j}(x_{i,j})e_{s,t}$%
\ where $\Delta _{s,t}^{i,j}$\ are additive maps from $I_{i,j}$\ to $I_{s,t}.
$

\bigskip

The following three lemmas are based on the results from
\cite{SayiKuzu}.

\bigskip 

\textbf{Lemma 1} Let $\Delta $\ be an arbitrary Lie derivation of $R.$\ Then%
\begin{eqnarray}
\Delta (xe_{i+1,i}) &=&\sum \Delta
_{i+1,t}^{i+1,i}(x)e_{i+1,t}+\sum_{s\neq i+1}\Delta _{s,i}^{i+1,i}(x)e_{s,i}
\\
&&+\sum_{k\neq i,i+1}\Delta _{k,k}^{i+1,i}(x)e_{k,k}+\Delta
_{n,1}^{i+1,i}(x)e_{n,1},  \notag \\
\Delta (ye_{1,n}) &=&\sum \Delta _{1,t}^{1,n}(y)e_{1,t}+\sum_{s\neq
1}\Delta _{s,n}^{1,n}(y)e_{s,n} \\
&&+\sum_{k\neq 1,n}\Delta _{k,k}^{1,n}(y)e_{k,k}+\Delta
_{n-1,1}^{1,n}(y)e_{n-1,1}  \notag \\
&&+\Delta _{n-1,2}^{1,n}(y)e_{n-1,2}+\Delta _{n,1}^{1,n}(y)e_{n,1}+\Delta
_{n,2}^{1,n}(y)e_{n,2}  \notag \\
&&  \notag \\
\text{for}\ 1 &<&i<n-1\ \text{and}  \notag \\
&&  \notag \\
\Delta _{k,k}^{i,j} &=&\Delta _{m,m}^{i,j}:I_{i,j}\rightarrow C(K)
\end{eqnarray}%
for $1\leq i,j\leq n$, \ $k\neq i,j$\ \ and \ $m\neq i,j.$

\bigskip 

\textbf{Proof} The equality $x_{i,j}e_{i,j}\ast y_{k,m}e_{k,m}=0$\ for $%
k\neq j$\ and $m\neq i$\ is equivalent to%
\begin{eqnarray*}
0 &=&\sum \Delta _{s,k}^{i,j}(x_{i,j})y_{k,m}e_{s,m}-\sum y_{k,m}\Delta
_{m,t}^{i,j}(x_{i,j})e_{k,t} \\
&&+\sum x_{i,j}\Delta _{j,t}^{k,m}(y_{k,m})e_{i,t}-\sum \Delta
_{s,i}^{k,m}(y_{k,m})x_{i,j}e_{s,j}
\end{eqnarray*}%
and (k,m) entry of this equality gives (3) for $k>m,$\ $\ k\neq i,j$\ and $%
m\neq i,j$. On the other hand, for $m\neq j,$\ $s\neq k,i$\ and $y_{k,m}=1,$%
\ we have $\Delta _{s,k}^{i,j}=0$. Similarly, for $k>m,$\ $\ k\neq i,$\ $\
t\neq m,j$\ and $y_{k,m}=1,$\ we get$\ \Delta _{m,t}^{i,j}=0.$\ Therefore,
we obtain (1) and (2). In particular, we get%
\begin{eqnarray*}
\Delta (xe_{2,1}) &=&\sum \Delta _{2,t}^{2,1}(x)e_{2,t}+\sum_{s\neq
2}\Delta _{s,1}^{2,1}(x)e_{s,1} \\
&&+\sum_{k\neq 1,2}\Delta _{k,k}^{2,1}(x)e_{k,k}+\Delta
_{n,2}^{2,1}(x)e_{n,2}+\Delta _{n,3}^{2,1}(x)e_{n,3} \\
\Delta (xe_{n,n-1}) &=&\sum \Delta _{n,t}^{n,n-1}(x)e_{n,t}+\sum_{s\neq
n}\Delta _{s,n-1}^{n,n-1}(x)e_{s,n-1} \\
&&+\sum_{k\neq n,n-1}\Delta _{k,k}^{n,n-1}(x)e_{k,k}+\Delta
_{n-2,1}^{n,n-1}(x)e_{n-2,1} \\
&&+\Delta _{n-1,1}^{n,n-1}(x)e_{n-1,1}.
\end{eqnarray*}

\textbf{Lemma 2} Let $\Delta :R\rightarrow R$\ be a Lie derivation. Then
there exist a diagonal derivation $\sigma _{d}$\ and two inner derivations $%
I_{A}$, $I_{B}$\ such that $(\Delta -\sigma _{d}-I_{A}-I_{B})(e_{i+1,i})$\
has nonzero entries only on i+1 row, main diagonal and (n,1) entry.

\bigskip 

\textbf{Proof} Let $d_{i+1}=\sum_{k=1}^{i}\Delta _{k+1,k}^{k+1,k}(1)$\ and $%
d=\sum_{i=2}^{n}d_{i}e_{i,i}$\ be a diagonal matrix. If $\sigma _{d}$\ is
the diagonal derivation induced by the matrix $d,$\ then%
\begin{eqnarray*}
\sigma _{d}(e_{i+1,i}) &=&de_{i+1,i}-e_{i+1,i}d \\
&=&\Delta _{i+1,i}^{i+1,i}(1)e_{i+1,i}.
\end{eqnarray*}%
So the $(i+1,i)$\ entry of the matrix $(\Delta -\sigma _{d})(e_{i+1,i})$\ is
zero. Define a matrix $A=\sum_{i=2}^{n-1}(-\Delta
_{i+1,1}^{i+1,i}(1))e_{i,1}$\ and a matrix$\ B=[b_{i,j}]$\ such that $%
b_{u,u}=0,$\ $b_{j,1}=0$, $b_{v,i+1}=\Delta _{v,i}^{i+1,i}(1)$\ for $i<n$, $%
v\neq i,i+1.$\ Let $I_{A}$\ and $I_{B}$\ be the inner derivations induced by the
matrices $A$\ and $B,$\ respectively. If we use $%
I_{A}(e_{i+1,i})=Ae_{i+1,i}-e_{i+1,i}A$\ and $%
I_{B}(e_{i+1,i})=Be_{i+1,i}-e_{i+1,i}B,$\ we obtain (i+1,1) entry and all
entries of $(\Psi -I_{A}-I_{B})(e_{i+1,i})$\ on i-th column are zeros. The
lemma is proved. In particular, for $\Psi =\Delta -\sigma _{d}-I_{A}-I_{B}$,
we have%
\begin{eqnarray*}
\Psi (e_{2,1}) &=&\sum_{t\neq 1}\Psi
_{i+1,t}^{2,1}(1)e_{i+1,t}+\sum_{k\neq 1}\Psi _{k,k}^{2,1}(1)e_{k,k} \\
&&+\Psi _{n,2}^{2,1}(1)e_{n,2}+\Psi _{n,3}^{2,1}(1)e_{n,3}\text{ \ ,} \\
\\
\Psi (e_{n,n-1}) &=&\sum_{t\neq 1,n-1}\Psi
_{n,t}^{n,n-1}(1)e_{n,t}+\sum_{k\neq n-1}\Psi _{k,k}^{n,n-1}(1)e_{k,k} \\
&&+\Psi _{n-1,1}^{n,n-1}(1)e_{n-1,1}+\Psi _{n-2,1}^{n,n-1}(1)e_{n-2,1}.
\end{eqnarray*}

\bigskip 

\bigskip 

By using the following relations%
\begin{eqnarray*}
\Psi (e_{2,1}\ast e_{i+1,i}) &=&0\text{ }(2<i<n) \\
\Psi (e_{i+1,i}\ast e_{j+1,j}) &=&0\text{ }(i\neq 1,\text{ \ }j\neq 1,\text{
\ }i\neq j-1,j,j+1) \\
\Psi (e_{n,n-1}\ast e_{i+1,i}) &=&0\text{ }(1<i<n-2) \\
\Psi (e_{i+1,i}\ast e_{i,i-1}) &=&\Psi (e_{i+1,i-1}) \\
\Psi (e_{i+1,i-1}\ast e_{i,i-1}) &=&0 \\
\Psi (e_{i+1,i-1}\ast e_{i+1,i}) &=&0
\end{eqnarray*}%
we get%
\begin{eqnarray*}
\Psi (e_{i+1,i}) &=&\sum_{k\neq i}\Psi _{k,k}^{i+1,i}(1)e_{k,k}+\Psi
_{n,1}^{i+1,i}(1)e_{n,1}\ (i\neq 1,\text{ }n-1), \\
\Psi (e_{2,1}) &=&\sum_{k\neq 1}\Psi _{k,k}^{2,1}(1)e_{k,k}+\Psi
_{n,2}^{2,1}(1)e_{n,2}, \\
\Psi (e_{n,n-1}) &=&\sum_{k\neq n-1}\Psi _{k,k}^{n,n-1}(1)e_{k,k}+\Psi
_{n-1,1}^{n,n-1}(1)e_{n-1,1}.
\end{eqnarray*}

\bigskip 

\textbf{Lemma 3} There exists a central Lie derivation $\Omega $\ of $R$\
such that (k,k) entries of the matrix $(\Delta -\sigma
_{d}-I_{A}-I_{B}-\Omega )(xe_{i+1,i})$\ and (m,m) entries of $(\Delta
-\sigma _{d}-I_{A}-I_{B}-\Omega )(ye_{1,n})$\ are equal to zero for 
$i<n,$\ $\ k\neq i,i+1,$\ $\ m\neq 1,n.$

\textbf{Proof} Let $\Psi =\Delta -\sigma _{d}-I_{A}-I_{B}$\ for brevity. If $%
1<i<n,$ then the equality $0=\Psi (xe_{i+1,i}\ast ye_{1,n})$\ gives $\Psi
_{n,1}^{i+1,i}(K)\subset Ann_{K}(J).$\ Besides, we have $\Psi
_{n,1}^{i+1,i}(J)=0$\ by the relations $\Psi (xye_{i+1,i})=\Psi
(xe_{i+1,i}\ast ye_{i,i})$\ and $\Psi (yxe_{i+1,i})=\Psi (ye_{i+1,i+1}\ast
xe_{i+1,i}).$\ In addition, we get $\Psi _{n,1}^{2,1}(J)=0$\ and$\ \Psi
_{n,1}^{2,1}(K)\subset Ann_{K}(J)$\ by using the relations $\Psi
(xye_{2,1})=\Psi (xe_{2,1}\ast ye_{1,1}),$\ $\Psi (yxe_{2,1})=\Psi
(ye_{2,2}\ast xe_{2,1}),$\ $\Psi (ye_{1,2}\ast xe_{2,1})=\Psi
(yxe_{1,1}-xye_{2,2})$\ and $\Psi (xe_{2,1}\ast ye_{n,n})=0.$\ Moreover, we
obtain $\Psi _{n,1}^{n,n-1}(K)\subset Ann_{K}(J)$\ and $\Psi
_{n,1}^{n,n-1}(J)=0$\ by the equalities $\Psi (xye_{n,n-1})=\Psi
(xe_{n,n-1}\ast ye_{n-1,n-1}),$\ $\Psi (yxe_{n,n-1})=\Psi (ye_{n,n}\ast
xe_{n,n-1}),$\ $\Psi (ye_{n-1,n}\ast xe_{n,n-1})=\Psi
(yxe_{n-1,n-1}-xye_{n,n})$\ and $\Psi (xe_{n,n-1}\ast ye_{1,1})=0.$\
Furthermore, If we consider the relations $\Psi (ye_{1,2}\ast ze_{2,n})=\Psi
(yze_{1,n}),$\ $\Psi (ye_{1,n}\ast ze_{2,n})=0$\ and $\Psi (ye_{1,n}\ast
ze_{1,2})=0,$\ we have $\Psi _{n,1}^{1,n}(J^{2})=0$\ and $\Psi
_{n,1}^{1,n}(J)\subset Ann_{J}(K).$\ Now let\ $\alpha _{i}=\Psi
_{k,k}^{i+1,i}$\ and $\alpha _{n}=\Psi _{m,m}^{1,n}$\ where $m\neq 1,n,$\ $\
k\neq i,i+1$\ and $i<n.$\ We know that $\alpha _{i}=\Psi
_{k,k}^{i+1,i}:K\rightarrow C(K)$\ $(i<n,$\ \ $k\neq i,i+1)$\ and $\alpha
_{n}=\Psi _{m,m}^{1,n}:J\rightarrow C(K)$\ $(m\neq 1,n)$\ by (3) in Lemma 1. 
$\Psi (xye_{i+1,i})=\Psi (xe_{i+1,i}\ast ye_{i,i})$\ gives $\alpha _{i}(J)=0$%
\ and we get $\alpha _{n}(J^{2})=0$\ by $\Psi (yze_{1,n})=\Psi (ye_{1,n}\ast
ze_{n,n}).$\ On the other hand, $\Psi (xye_{i+1,i+1}-yxe_{i,i})=\Psi
(xe_{i+1,i}\ast ye_{i,i+1})$\ gives $\Psi _{n,1}^{i+1,i+1}=\Psi _{n,1}^{i,i}$%
\ and $\Psi _{n,1}^{i,i}(xy)=\Psi _{n,1}^{i,i}(yx)$\ for $1<i<n.$\ Say $%
u=\Psi _{n,1}^{i,i}$\ for $1<i<n.$\ Also we obtain $u:J\rightarrow Ann_{K}(J)
$ by the relation $ye_{2,2}\ast ze_{1,n}=0.$\ Then the map%
\begin{equation*}
\begin{array}{llll}
\Omega : & R & \rightarrow  & R \\ 
& [a_{i,j}] & \rightarrow  & 
\begin{array}{l}
\left( \sum_{i=1}^{n-1}\alpha _{i}(a_{i+1,i})+\alpha _{n}(a_{1,n})\right) .E
\\ 
+\left( 
\begin{array}{c}
\sum_{i=1}^{n}u(a_{i,i})+\sum_{i=1}^{n-1}\Psi _{n,1}^{i+1,i}(a_{i+1,i}) \\ 
+\Psi _{n,1}^{1,n}(a_{1,n})%
\end{array}%
\right) e_{n,1}%
\end{array}%
\end{array}%
\end{equation*}%
is the required central Lie derivation. The lemma is proved.

\bigskip 

\textbf{Lemma 4}\ There exists a ring derivation $\bar{\Pi}$\ of $R$\ such
that (i,j) entry of $(\Delta -\sigma _{d}-I_{A}-I_{B}-\Omega -\bar{\Pi}%
)(x_{i,j}e_{i,j})$\ is zero for all $i,j.$

\bigskip 

\textbf{Proof} Let $\Psi =\Delta -\sigma _{d}-I_{A}-I_{B}-\Omega $\ $\ $for
brevity$.$\ Then\ (i,k) entry of the relation $\Psi
(x_{i,j}y_{j,k}e_{i,k})=\Psi (x_{i,j}e_{i,j}\ast y_{j,k}e_{j,k})$\ gives $%
\Psi _{i,k}^{i,k}(x_{i,j}y_{j,k})=\Psi
_{i,j}^{i,j}(x_{i,j})y_{j,k}+x_{i,j}\Psi _{j,k}^{j,k}(y_{j,k})$\ for all $%
i,j,k$\ where $i\neq k.$\ This means all $\Psi _{u,v}^{u,v}$\ are equal to
each other and $\Psi _{u,v}^{u,v}$\ is a derivation of both $K$\ and $J.$\
Let $\Pi =\Psi _{u,v}^{u,v}.$\ Then the map%
\begin{equation*}
\begin{array}{llll}
\bar{\Pi}: & R & \rightarrow  & R \\ 
& [a_{i,j}] & \rightarrow  & [\Pi (a_{i,j})]%
\end{array}%
\end{equation*}%
is the ring derivation we need. The lemma is proved.

\bigskip 

Now we have the following equalities for $\Psi =\Delta -\sigma
_{d}-I_{A}-I_{B}-\Omega -\bar{\Pi}.$%
\begin{eqnarray*}
\Psi (e_{2,1}) &=&\Psi _{1,1}^{2,1}(1)e_{1,1}+\Psi
_{2,2}^{2,1}(1)e_{2,2}+\Psi _{n,2}^{2,1}(1)e_{n,2}, \\
\Psi (e_{i+1,i}) &=&\Psi _{i,i}^{i+1,i}(1)e_{i,i}+\Psi
_{i+1,i+1}^{i+1,i}(1)e_{i+1,i+1}\text{ \ }(1<i<n-1), \\
\Psi (e_{n,n-1}) &=&\Psi _{n-1,n-1}^{n,n-1}(1)e_{n-1,n-1}+\Psi
_{n,n}^{n,n-1}(1)e_{n,n}+\Psi _{n-1,1}^{n,n-1}(1)e_{n-1,1}.
\end{eqnarray*}

\bigskip 

\textbf{Lemma 5} There exists a special Lie derivation $\bar{\Delta}$\ of $R$%
\ such that (n,2) entry of $(\Delta -\sigma _{d}-I_{A}-I_{B}-\Omega -\bar{\Pi%
}-\bar{\Delta})(xe_{2,1})$\ and (n-1,1) entry of $\ (\Delta -\sigma
_{d}-I_{A}-I_{B}-\Omega -\bar{\Pi}-\bar{\Delta})(xe_{n,n-1})$\ are zeros$.$

\bigskip 

\textbf{Proof} Let $\Psi =\Delta -\sigma _{d}-I_{A}-I_{B}-\Omega -\bar{\Pi}$%
\ for brevity$.$\ The relations $\Psi (x_{1}e_{2,1}\ast x_{2}e_{2,1})=0$\
and $\Psi (x_{1}e_{n,n-1}\ast x_{2}e_{n,n-1})=0$\ gives $\Psi
_{n,2}^{2,1}(x_{1})x_{2}-\Psi _{n,2}^{2,1}(x_{2})x_{1}=0$\ and $x_{1}\Psi
_{n-1,1}^{n,n-1}(x_{2})-x_{2}\Psi _{n-1,1}^{n,n-1}(x_{1})=0,$\ respectively.
Besides, the relations $\Psi (xye_{2,1})=\Psi (xe_{2,1}\ast ye_{1,1})$, $%
\Psi (yxe_{2,1})=\Psi (ye_{2,2}\ast xe_{2,1})$\ and $\Psi (xe_{2,1}\ast
ye_{n,n})=0$\ gives $\Psi _{n,2}^{2,1}(J)=0$\ and$\ \Psi
_{n,2}^{2,1}(K)\subset Ann_{K}(J)$. Moreover, we obtain $\Psi
_{n-1,1}^{n,n-1}(J)=0$\ and $\Psi _{n-1,1}^{n,n-1}(K)\subset Ann_{K}(J)$\ by
the relations $\Psi (yxe_{n,n-1})=\Psi (ye_{n,n}\ast xe_{n,n-1}),$\ $\Psi
(xe_{n,n-1}\ast ye_{1,1})=0$\ and $\Psi (xye_{n,n-1})=\Psi (xe_{n,n-1}\ast
ye_{n-1,n-1}),$\ respectively$.$\ Then the map%
\begin{equation*}
\begin{array}{llll}
\bar{\Delta}: & R & \rightarrow  & R \\ 
& [a_{i,j}] & \rightarrow  & \Psi _{n,2}^{2,1}(a_{2,1})e_{n,2}+\Psi
_{n-1,1}^{n,n-1}(a_{n,n-1})e_{n-1,1}%
\end{array}%
\end{equation*}%
is a special Lie derivation of type I defined in (B1). The lemma is proved.

\bigskip 

\textbf{Lemma 6} There exists an inner derivation $I_{C}$\ such that $%
(\Delta -\sigma _{d}-I_{A}-I_{B}-\Omega -\bar{\Pi}-\bar{\Delta}%
-I_{C})(e_{i+1,i})$\ has zero (i,i) entry for all $i<n.$

\bigskip 

\textbf{Proof} Let $\Psi =\Delta -\sigma _{d}-I_{A}-I_{B}-\Omega -\bar{\Pi}-%
\bar{\Delta}$\ for brevity and$\ C=\sum_{k=1}^{n-1}\Psi
_{k,k}^{k+1,k}(1)e_{k,k+1}.$\ Then the inner derivation $I_{C}$\ induced by
the matrix $C$\ is the required inner derivation. Note that the positions
with nonzero entries of $I_{C}(e_{i,j})$\ and $(\Delta -\sigma
_{d}-I_{A}-I_{B}-\Omega -\bar{\Pi}-\bar{\Delta})(e_{i,j})$\ are the same for 
$i>j.$

\bigskip 

\bigskip 

Now let $\Psi =\Delta -\sigma _{d}-I_{A}-I_{B}-\Omega -\bar{\Pi}-\bar{\Delta}%
-I_{C}$\ for brevity. Then $\Psi (e_{i+1,i})=\Psi
_{i+1,i+1}^{i+1,i}(1)e_{i+1,i+1}$\ for $i<n.$\ But we have $\Psi
(e_{i+1,i-1})=0$\ and then $\Psi (e_{i+1,i})=0$\ by the relations $\Psi
(e_{i+1,i}\ast e_{i,i-1})=\Psi (e_{i+1,i-1})$\ and $\Psi (e_{i+1,i-1}\ast
e_{i,i-1})=0$. It means $\Psi (e_{i,j})=0$\ for all $i>j.$

\bigskip 

For any numbers $k$\ and $m$\ such that $k>m,$\ $\ k\neq 1,i$, $\ m\neq
i+1,n,$\ we obtain%
\begin{eqnarray*}
\Psi (xe_{i+1,i}) &=&\Psi _{i+1,1}^{i+1,i}(x)e_{i+1,1}+\Psi
_{n,i}^{i+1,1}(x)e_{n,i}\text{ \ }(1<i<n-1) \\
\Psi (xe_{2,1}) &=&\Psi _{n,3}^{2,1}(x)e_{n,3} \\
\Psi (xe_{n,n-1}) &=&\Psi _{n-2,1}^{n,n-1}(x)e_{n-2,1}
\end{eqnarray*}%
by the relation $\Psi (xe_{i+1,i}\ast e_{k,m})=\Psi (xe_{i+1,i})\ast
e_{k,m}=0.\ $Besides, for$\ $\ $2<i<n-1,$\ the relation $\Psi
(x_{1}x_{2}e_{i+1,i-1})=\Psi (x_{1}e_{i+1,i}\ast x_{2}e_{i,i-1})$\ gives $%
\Psi _{n,i}^{i+1,i}=0$\ and $\Psi _{i,1}^{i,i-1}=0.$\ In addition, we get $%
\Psi _{n,2}^{3,2}=0=\Psi _{n,3}^{2,1}$\ and $\Psi _{n-1,1}^{n-1,n-2}=0=\Psi
_{n-2,1}^{n,n-1}$\ by the equalities $\Psi (x_{1}x_{2}e_{3,1})=\Psi
(x_{1}e_{3,2}\ast x_{2}e_{2,1})$\ and $\Psi (x_{1}x_{2}e_{n,n-2})=\Psi
(x_{1}e_{n,n-1}\ast x_{2}e_{n-1,n-2})$. Finally, we have%
\begin{equation}
\Psi (xe_{i+1,i})=0\text{ and so }\Psi (xe_{k,m})=0  \tag{$\ast $}
\end{equation}%
for any $i<n$\ and $k>m.$\ Also we obtain%
\begin{eqnarray*}
\Psi (ye_{1,n}) &=&\Psi _{1,1}^{1,n}(y)e_{1,1}+\Psi
_{n-1,1}^{1,n}(y)e_{n-1,1}+\Psi _{n-1,2}^{1,n}(y)e_{n-1,2} \\
&&+\Psi _{n,2}^{1,n}(y)e_{n,2}+\Psi _{n,n}^{1,n}(y)e_{n,n}
\end{eqnarray*}%
by the relations $\Psi (ye_{1,n-1})=\Psi (ye_{1,n}\ast e_{n,n-1})$\ and $%
\Psi (ye_{2,n})=\Psi (e_{2,1}\ast ye_{1,n}).$

\bigskip 

\textbf{Lemma 7}\ There exists an almost annihilator derivation $\Gamma $\
of $R$\ such that (n,i) entry of $(\Delta -\sigma _{d}-I_{A}-I_{B}-\Omega -%
\bar{\Pi}-\bar{\Delta}-I_{C}-\Gamma )(ye_{1,i})$\ and (i,1) entry of $%
(\Delta -\sigma _{d}-I_{A}-I_{B}-\Omega -\bar{\Pi}-\bar{\Delta}-I_{C}-\Gamma
)(ye_{i,n})$\ are zeros for $1\leq i\leq n$.

\bigskip 

\textbf{Proof}\ Let $\Psi =\Delta -\sigma _{d}-I_{A}-I_{B}-\Omega -\bar{\Pi}-%
\bar{\Delta}-I_{C}$\ for brevity. Then we obtain $\Psi _{n,n}^{1,n}=\Psi
_{n,i}^{1,i}$, $\Psi _{n,n}^{1,n}(yx)=\Psi _{n,n}^{1,n}(y)x$\ and $\Psi
_{1,1}^{1,n}=\Psi _{i,1}^{i,n},$\ $\Psi _{1,1}^{1,n}(xy)=x\Psi
_{1,1}^{1,n}(y)$\ by the relations $\Psi (ye_{1,n}\ast xe_{n,i})=\Psi
(yxe_{1,i})$\ and $\Psi (xe_{i,1}\ast ye_{1,n})=\Psi (xye_{i,n})$\ for $%
1<i<n.$\ Furthermore, we get $\Psi _{n,1}^{1,1}=\Psi _{n,2}^{1,2}$\ and $%
\Psi _{n,1}^{n,n}=\Psi _{n-1,1}^{n-1,n}$\ by $\Psi (xe_{2,1}\ast
ye_{1,2})=\Psi (xye_{2,2}-yxe_{1,1})$\ and $\Psi (xe_{n,n-1}\ast
ye_{n-1,n})=\Psi (xye_{n,n}-yxe_{n-1,n-1}).$\ Now we have $\Psi
_{n,n}^{1,n}=\Psi _{n,i}^{1,i}$\ and $\Psi _{1,1}^{1,n}=\Psi _{i,1}^{i,n}$.
Let $\lambda =\Psi _{1,1}^{1,n}=\Psi _{i,1}^{i,n}$\ and $\mu =\Psi
_{n,n}^{1,n}=\Psi _{n,i}^{1,i}.$\ Now consider the relation $\Psi
(yze_{1,n})=\Psi (ye_{1,2}\ast ze_{2,n}).$\ This relation gives $0=z\mu
(y)+\lambda (z)y.$\ Then the map%
\begin{equation*}
\begin{array}{llll}
\Gamma : & R & \rightarrow  & R \\ 
& ye_{1,n} & \rightarrow  & \lambda (y)e_{1,1}+\mu (y)e_{n,n} \\ 
& ye_{i,n} & \rightarrow  & \lambda (y)e_{i,1}\text{ \ }(i>1) \\ 
& ye_{1,j} & \rightarrow  & \mu (y)e_{n,j}\text{ \ }(j<n) \\ 
& x_{i,j}e_{i,j} & \rightarrow  & 0\text{ \ }((i,j)\neq (1,k),(m,n)\text{
for any }k,m)%
\end{array}%
\end{equation*}%
is an almost annihilator derivation of $R$\ and this proves the lemma.

\bigskip 

\textbf{Lemma 8}\ There exists a special Lie derivation $\Theta $\ of $R$\
such that $\Delta -\sigma _{d}-I_{A}-I_{B}-\Omega -\bar{\Pi}-\bar{\Delta}%
-I_{C}-\Gamma -\Theta $\ is equal to zero map.

\bigskip 

\textbf{Proof}\ Let $\Psi =\Delta -\sigma _{d}-I_{A}-I_{B}-\Omega -\bar{\Pi}-%
\bar{\Delta}-I_{C}-\Gamma $\ for brevity$.$\ Then the relations $\Psi
(yxe_{1,n-1})=\Psi (ye_{1,n}\ast xe_{n,n-1}),$\ $\Psi (xye_{2,n})=\Psi
(xe_{2,1}\ast ye_{1,n})$\ and $\Psi (yxe_{2,n-1})=\Psi (ye_{2,n}\ast
xe_{n,n-1})$\ give $\Psi _{n,1}^{1,n-1}(yx)=-x\Psi _{n-1,1}^{1,n}(y),$\ $%
\Psi _{n,2}^{1,n-1}(yx)=-x\Psi _{n-1,2}^{1,n},$\ $\Psi
_{n-1,1}^{2,n}(xy)=-\Psi _{n-1,2}^{1,n}(y)x,$\ $\Psi _{n,1}^{2,n}(xy)=-\Psi
_{n,2}^{1,n}(y)x$\ and$\ \Psi _{n,1}^{2,n-1}(yx)=-x\Psi _{n-1,1}^{2,n}(y).$\
Let $\alpha =\Psi _{n-1,1}^{1,n},$\ $\beta =\Psi _{n-1,2}^{1,n}$\ and $%
\gamma =\Psi _{n,2}^{1,n}.$\ Then one can see that $\alpha (yx)=x\alpha (y),$%
\ $\beta (yx)=x\beta (y),$\ $\beta (xy)=\beta (y)x,$\ $\gamma (xy)=\gamma
(y)x,$\ $\Psi _{n,1}^{1,n-1}=-\alpha,$\ $\Psi_{n,2}^{1,n-1}=\Psi_{n-1,1}^{2,n}=-\Psi_{n,1}^{2,n-1}=-\beta,$\ $\Psi_{n,1}^{2,n}=-\gamma.$

If $n>4$\ and $2<i<n-1,$\ then we obtain $\alpha (J^{2})=\beta
(J^{2})=\gamma (J^{2})=0\ $and that$\ \alpha (y),\beta (y),\gamma (y)\subset
Ann_{K}(J)$\ by the relations $\Psi (yze_{1,n})=\Psi (ye_{1,i}\ast ze_{i,n}),
$\ $\Psi (yze_{1,n})=\Psi (ye_{1,2}\ast ze_{2,n}),$\ $\Psi (yze_{1,n})=\Psi
(ye_{1,n-1}\ast ze_{n-1,n}),$\ $\Psi (yze_{1,n})=\Psi (ye_{1,1}\ast ze_{1,n})
$\ and $0=\Psi (ye_{1,n}\ast ze_{2,n})$.

If $n=4,$\ then the relations $\Psi (yze_{1,4})=\Psi (ye_{1,j}\ast ze_{j,4})$%
\ and $\Psi (ye_{1,4}\ast ze_{2,4})=0$\ gives $\alpha (J^{2})=\beta
(J^{2})=\gamma (J^{2})=0\ $and$\ \alpha (y),\beta (y),\gamma (y)\subset
Ann_{K}(J)$\ for $1\leq j\leq 4.$\ So $\alpha ,\beta $\ and $\gamma $\
satisfies the conditions given in (B2) and the map%
\begin{equation*}
\begin{array}{llll}
\Theta : & R & \rightarrow  & R \\ 
& ye_{1,n} & \rightarrow  & \alpha (y)e_{n-1,1}+\beta (y)e_{n-1,2}+\gamma
(y)e_{n,2} \\ 
& ye_{1,n-1} & \rightarrow  & -\alpha (y)e_{n,1}-\beta (y)e_{n,2} \\ 
& ye_{2,n} & \rightarrow  & -\beta (y)e_{n-1,1}-\gamma (y)e_{n,1} \\ 
& ye_{2,n-1} & \rightarrow  & \beta (y)e_{n,1} \\ 
& x_{i,j}e_{i,j} & \rightarrow  & 0\text{ \ \ \ }((i,j)\neq
(1,n-1),(1,n),(2,n),(2,n-1))%
\end{array}%
\end{equation*}%
becomes a special Lie derivation of type II$.$\ For the last part, we need to
show $(\Delta -\sigma _{d}-I_{A}-I_{B}-\Omega -\bar{\Pi}-\bar{\Delta}%
-I_{C}-\Gamma -\Theta )(x_{i,j}e_{i,j})=0$\ for all i,j. Let $\Xi =\Delta
-\sigma _{d}-I_{A}-I_{B}-\Omega -\bar{\Pi}-\bar{\Delta}-I_{C}-\Gamma -\Theta 
$\ for brevity. Then $\Xi (ye_{i,j})=0$\ for $i<j$\ since we have%
\begin{eqnarray*}
\Xi (e_{i,1}) &=&0, \\
\Xi (ye_{1,j}) &=&0, \\
\Xi (ye_{i,j}) &=&\Xi (e_{i,1}\ast ye_{1,j}) \\
&=&\Xi (e_{i,1})\ast ye_{1,j}+e_{i,1}\ast \Xi (ye_{1,j}).
\end{eqnarray*}%
Besides, we know that $\Xi (xe_{i,j})=0$\ for $i>j$\ by $(\ast ).$\ It
suffices to show that $\Xi (ye_{i,i})=0$\ for all i. If $i>j$, then $\Xi
(ye_{i,i}-ye_{j,j})=\Xi (e_{i,j}\ast ye_{j,i})=0$\ as $\Xi (ye_{k,m})=0$\
for $k<m.$\ This means $\Xi (ye_{i,i})=\Xi (ye_{j,j})$\ for all i,j. So it
is enough to check $\Xi (ye_{1,1})$\ is equal to $0.$\ Consider the products 
$\Xi (ye_{1,1}\ast e_{i,j})=\Xi (ye_{1,1})\ast e_{i,j}=0$\ for $i>j.$\ These
products give that entries on the $i-th$\ column and $j-th$\ row of $\Xi
(ye_{1,1})$\ are zeros for $1<j<n$\ and $i>2$\ which means $\Xi
(ye_{1,1})=\Xi _{1,2}^{1,1}(y)e_{1,2}+\Xi _{n,2}^{1,1}(y)e_{n,2}.$\ Finally,
the relation $0=\Xi (ye_{2,1})=\Xi (e_{2,1}\ast ye_{1,1})=e_{2,1}\ast \Xi
(ye_{1,1})$\ gives $\Xi _{1,2}^{1,1}=\Xi _{n,2}^{1,1}=0.$\ The lemma is
proved.

\bigskip 

The proof of the main theorem follows by Lemma 1-8 for $n\geq 4$.

\bigskip 

\textbf{Discussion for n=3:}

\bigskip 

Let $K$\ be commutative now$.$\ Then Lemma 1-5 applies. Moreover, the
relations $yze_{1,3}=ye_{1,1}\ast ze_{1,3}$, $yze_{1,3}=ye_{1,2}\ast
ze_{2,3},$\ \ $yze_{1,3}=ye_{1,3}\ast ze_{3,3},$\ \ $xye_{2,3}=xe_{2,1}\ast
ye_{1,3},$\ \ $yxe_{1,2}=ye_{1,3}\ast xe_{3,2}$\ gives all the conditions of
special Lie derivation $\Lambda $\ defined in (B3). Let $\Psi =\Delta
-\sigma _{d}-I_{A}-I_{B}-\Omega -\bar{\Pi}-\bar{\Delta}-\Lambda $\ for
brevity. Then we have $\Psi (xe_{2,1})=\Psi _{1,1}^{2,1}(x)e_{1,1}+\Psi
_{2,2}^{2,1}(x)e_{2,2},$\ \ $\Psi (xe_{3,2})=\Psi
_{2,2}^{3,2}(x)e_{2,2}+\Psi _{3,3}^{3,2}(x)e_{3,3}$\ and$\ \Psi
(xe_{3,1})=\Psi _{2,1}^{3,1}(x)e_{2,1}+\Psi _{3,2}^{3,1}(x)e_{3,2}.$\  Now
Lemma 6 applies and then we find that the images of $xe_{i+1,i}$\ under $%
\Delta -\sigma _{d}-I_{A}-I_{B}-\Omega -\bar{\Pi}-\bar{\Delta}-\Lambda -I_{C}
$\ \ are zeros by the equalities $x_{1}x_{2}e_{3,1}=x_{1}e_{3,2}\ast
x_{2}e_{2,1},$\ \ $x_{1}e_{2,1}\ast x_{2}e_{3,1}=0$\ and $x_{1}e_{3,2}\ast
x_{2}e_{3,1}=0.$\ After applying Lemma 7, we also obtain the image of $%
ye_{1,3}$\ is zero under $\Delta -\sigma _{d}-I_{A}-I_{B}-\Omega -\bar{\Pi}-%
\bar{\Delta}-\Lambda -I_{C}-\Gamma .$

Now let $\Xi =\Delta -\sigma _{d}-I_{A}-I_{B}-\Omega -\bar{\Pi}-\bar{\Delta}%
-\Lambda -I_{C}-\Gamma $\ and consider the relations $ye_{1,2}=ye_{1,3}\ast
e_{3,2},$\ \ $ye_{2,3}=e_{2,1}\ast ye_{1,3}$. Then we find $\Xi (ye_{i,j})=0$%
\ for $i<j.$\ Furthermore, we obtain $\Xi (ye_{k,k})=0$\ by $e_{i,j}\ast
ye_{k,k}=0.$\ So our main theorem holds for n=3 if $K$\ is commutative. Note
that the ring $K$\ is chosen to be commutative to obtain central Lie
derivations of $R_{3}(K,J).$


\end{document}